\documentclass[a4paper,12pt]{article}
\usepackage[cp1251]{inputenc}
\usepackage[english]{babel}
\usepackage{amsmath}
\usepackage{mathrsfs}
\usepackage{amsfonts}
\usepackage[all]{xy}
\usepackage{amscd}
\usepackage{amssymb}
\begin{document}

\title{On geometric properties of the functors of positively homogenous and semiadditive functionals}

\author{Lesya Karchevska}
\date{}
\maketitle

Department of Mechanics and Mathematics, Ivan Franko National University of Lviv,
Universytetska st.,1 79602 Lviv Ukraine, e-mail:
\emph{crazymaths@ukr.net}.

\begin{abstract}
In this paper we investigate the functors of $OH$ of positively
homogenous functionals and $OS$ of semiadditive functionals. We show
that $OH(X)\in AR$ if and only if $X$ is openly generated, and
$OS(X)\in AR$ if and only if $X$ is an openly generated compactum of
weight $\le \omega_1$. In section 3 we investigate the
multiplication maps of monads generated by the abovementioned
functors and consider when these mappings are soft.
\end{abstract}



\vspace*{\parindent} \vspace*{\parindent}

{\bf 0. Introduction.} V.Fedorchuk posed a general problem
concerning geometric properties of functors, that is, how functors
affect certain geometric properties of spaces and mappings between
them [11]. Under geometric properties we understand the property of
being an $AR$ for a space, the properties of being soft or a
Tychonov fibering for a mapping etc.

There were many investigations in this direction involving such
functors as the hyperspace functor $\exp$, the probability measures
functor $P$, the superextension functor $\lambda$, the inclusion
hyperspace functor $G$ and others (see, e.g. [10] or [11]).

Let us now consider as an example the functors of probability
measures $P$ and superextension $\lambda$. There is a natural
structure of linear convexity on $P(X)$. As for $\lambda$, de Groot
constructed some abstract convexity (not linear) on any space of the
form $\lambda (X)$ (see [12]), and this convexity is binary, whereas
the linear convexity on $P(X)$ is not.

The functors $\lambda$ and $P$ differ in their geometric properties
as well. Consider the property of being an $AR$, for instance. In
the metrizable case, $\lambda (X)\in AR$ if and only if $X$ is a
continuum, and $P(X)$ is an absolute retract for each compactum $X$.
When $X$ is not metrizable, the space $P(X)$ can be $AR$ only in
case $X$ is openly generated and of weight $\le \omega_1$. As for
the superextension functor, $\lambda(X)\in AR$ whenever $X$ is an
openly generated continuum, without limitations on weight.

The algebraic aspects of functors are formalized by the notion of a
monad in the sense of Eilenberg and Moore [13].

The notion of convexity considered in this paper is considerably
broader than the classic one: specifically, it is not restricted to
the context of linear spaces. Such convexities appeared in the
process of studying different structures like partially ordered
sets, semilattices, lattices, superextensions etc. We base our
approach on the notion of topological convexity from [14] where the
general convexity theory is covered from axioms to application in
different areas. T.Radul assigned to each monad $\mathbb{F}$ some
abstract convexity structure on every space $FX$, where $F$ is the
functorial part of the monad $\mathbb{F}$. Some additional
conditions on these monads (that they are $L$-monads which weakly
preserve preimages) guarantee that the considered convexities
generate the topology of the space $FX$ for the functor $F$ included
in an $L$-monad. It was shown that $L$-monads which weakly preserve
preimages and with binary convexities can give absolute retracts in
all weights [3]. Also, the morphisms of their algebras can be soft
in nonmetrizable case under certain conditions. Note that the
property of binarity of the convexity generated by monad
$\mathbb{F}$ is equivalent to the superextension monad being the
submonad of $\mathbb{F}$ (again [3]).

In this article we consider functors $OS$ and $OH$ (introduced in
[5], [6]), which both generate $L$-monads. The monad $\mathbb{OS}$
does not generate binary convexities, in turn $\mathbb{OH}$ does,
and this as well appears to be the reason for the difference in
their geometric properties: the properties of $OS$ are close to that
of $P$, and $OH$ is closer to $\lambda$.

\vspace*{\parindent}

{\bf 1. Definitions and facts.} In the present paper we shall deal
with objects and morphisms of the category $Comp$, that is, with
compact Hausdorff spaces and continuous mappings.

By $C(X)$, where $X\in Comp$, we denote the Banach space of all
continuous real-valued functions on $X$ with the sup-norm
$\Vert\varphi\Vert=\sup\lbrace\vert\varphi(x)\vert \ \vert \ x\in
X\rbrace$. By $c_X$, where $c\in \mathbb{R}$, we denote the constant
function: $c_X(x)=c$ for all $x\in X$.

Let $X\subset Y$. We say that a space $X$ is a \emph{retract} of $Y$
if there exists a map $r:Y\to X$ such that $r\vert_X =
\mathrm{id}_X$. The space $X$ is an \emph{absolute retract} (shortly
$X\in AR$), if for any embedding $i:X\hookrightarrow Y$ the subspace
$i(X)$ is a retract of $Y$.

Recall that a \emph{$\tau$- system}, where $\tau $ is any cardinal
number, is a continuous inverse system consisting of compacta of
weight $\le \tau$ and epimorphisms over a $\tau$-complete indexing
set. As usual, $\omega$ stands for the countable cardinal number. A
compactum $X$ is called \emph{openly generated}, if it can be
represented as the limit of some $\omega$-system with open bonding
mappings [1].

The mapping $f:X\to Y$ is called \emph{soft} if for any space $Z$
and its closed subset $A$, any functions $\psi:A\to X$, $\Psi:Z\to
Y$ with $\Psi\vert A = f\circ \psi$ there is a mapping $G:Z\to X$
such that $G\vert_A = \psi$ and $\Psi=f\circ G$ [1].

We say that a commutative diagram
\begin{displaymath}
\xymatrix{
  X \ar[r]^{p} \ar[d]_{q} & Y
   \ar[d]^{f}\\
  Z \ar[r]^{g} & T}
\end{displaymath}
is soft, if its characteristic map $\chi:X\to Y\times_T Z =
\{(y,z)\in Y\times Z \ \vert f(y)=g(z)\}$ defined by
$\chi(x)=(p(x),q(x))$ is soft.

A triple $(F,\eta,\mu)$, where $F$ is an endofunctor in category
$Comp$, $\eta: Id_{Comp}\to F$ and $\mu:F^2\to F$ are natural
transformations, is called \emph{monad} (in sense of Eilenberg and
Moore), if 1) $\mu X\circ \eta F(X) = \mu X\circ F(\eta X) =
\mathrm{id}_{F(X)}$; 2) $\mu X\circ \mu F(X) = \mu X\circ F(\mu X)$
[13].

Suppose that $\mathbb{F}=(F,\eta,\mu)$ is a monad. A pair $(X,\xi)$,
where $\xi:F(X)\to X$, is called an $\mathbb{F}$-algebra, if
$\xi\circ \eta X=id_X$ and $\xi\circ \mu X=\xi\circ F(\xi)$.

Let $\nu :C(X)\to \mathbb{R}$ be a functional. We say that $\nu$ is:
1) normed, if $\nu(1_X)=1$; 2) weakly additive, if for any $\phi \in
C(X)$ and $c\in \mathbb{R}$ we have $\nu(\phi+c_X)=\nu(\phi)+c$; 3)
order-preserving, whenever for any $\varphi,\psi\in C(X)$ such that
$\varphi(x)\le\psi(x)$ for all $x\in X$ (i.e. $\varphi\le \psi$) the
inequality $\nu(\varphi)\le\nu(\psi)$ holds; 4) positively
homogeneous, if for any $\varphi\in C(X)$ and any real $t\ge 0$ we
have $\nu(t\varphi)=t\nu(\varphi)$; 5) semiadditive, if
$\nu(\varphi+\psi)\le\nu(\varphi)+\nu(\psi)$.

Now for any space $X$ denote $V(X) = \prod_{\varphi\in
C(X)}[\min\varphi,\max\varphi]$. For any mapping $f:X\to Y$ let
$V(f)$ be a mapping such that $V(f)(\nu)(\varphi) = \nu(\varphi\circ
f)$ for any $\nu\in V(X), \ \varphi\in C(Y)$. Defined in that way,
$V$ forms a covariant functor in the category $Comp$.

For any space $X$ by $O(X)$ denote the set of functionals satisfying
1)--3) (order-preserving functionals), by $OH(X)$ the set of all
functionals on $C(X)$ which satisfy properties 1)--4) (positively
homogenous functionals), and by $OS(X)$ we denote the set of
functionals on $C(X)$ which satisfy properties 1)--5) (semiadditive
functionals).  Also recall that $P(X)$ stands for the set of all
functionals on $C(X)$ which are normed ($\Vert \mu \Vert = 1$),
positive ($\mu(\varphi)\ge 0$ for all $\varphi\ge 0$) and linear.
Let $F$ stand for one of $O,OH,OS,P$. The space $F(X)$ is considered
as the subspace of $V(X)$. For any function $f:X\to Y$, the map
$F(f):F(X)\to F(Y)$ is the restriction of $V(f)$ on the
corresponding space $F(X)$. Then $F$ forms a covariant functor in
$Comp$, which is a subfunctor of $V$.

It was shown in [5] and [6] that the functor $OS$ is normal, and
$OH$ is weakly normal, both $OH(X)$ and $OS(X)$ being convex
compacta for any space $X$.

Each of the abovementioned functors generates a monad. If $F$ is one
of $V,O,OH,OS,P$, the identity and multiplication maps are defined
as follows. The natural transformation $\eta:\mathrm{Id}_{{\it
Comp}}\to F$ is given by $\eta X(x)(\varphi)=\varphi(x)$ for any
$x\in X$ and $\varphi\in C(X)$, and the natural transformation $\mu:
F^2\to F$ given by $\mu X(\nu)(\varphi)=\nu(\pi_{\varphi})$, where
$\pi_{\varphi}:F(X)\to \mathbb{R}$,
$\pi_{\varphi}(\lambda)=\lambda(\varphi)$. Later by $\mu_F X$ we
shall denote the multiplication map for the corresponding functor
$F$. According to the characterization given in [15], by $L$-monad
we mean any submonad of $\mathbb{V}$. Hence, $\mathbb{OH}$ and
$\mathbb{OS}$, being submonads of $\mathbb{V}$ are both $L$-monads.

We say that an $L$-monad $\mathbb{F}=(F,\eta,\mu)$ {\it weakly
preserves preimages} ([3]) if for any mapping $f:X\to Y$ and any
closed subset $A\subset Y$ we have $\nu(\varphi)\in
[\min\varphi(f^{-1}(A)),\max\varphi(f^{-1}(A))]$ for all $\nu\in
(Ff)^{-1}(F(A))$ and all $\varphi\in C(X)$.

Let us recall the notion of convexities introduced in [3]. Let
$(F,\eta,\mu)$ be a monad, and $(X,\xi)$ be an $\mathbb{F}$-algebra.
Let $A$ be a closed subset of $X$. By $f_A$ denote the quotient map
$f_A:X\to X/A$, $a=f_A (A)$. We say that ${\cal {C}}_{\mathbb{F}}
(A)=\xi((Ff_{A}^{-1}(\eta(X/A)(a)))$ is the
$\mathbb{F}$-\emph{convex hull} of $A$. Also put ${\cal
C}_{\mathbb{F}}(\O)=\O$. The set $A$ is called
$\mathbb{F}$-\emph{convex} if $A = {\cal C}_{\mathbb{F}} (A)$.
Define ${\cal C}_{\mathbb F}(X,\xi)=\{A\subset X \ \vert \ A$ is
closed and $A={\cal C}_{\mathbb{F}}(A)\}$. The family ${\cal
C}_{\mathbb F}(X,\xi)$ forms a convexity on $X$. Also, any
$\mathbb{F}$-algebras morphism preserves convexities defined above
[3]. Later we'll restrict ourselves with the binary monads. A monad
$\mathbb{F}$ is \emph{binary} if
$\cal{C}_{\mathbb{F}}(\emph{X},\xi)$ is binary, i.e. the
intersection of each linked subsystem of
$\cal{C}_{\mathbb{F}}(\emph{X},\xi)$ is not empty (we call a family
of subsets of a space \emph{linked} if the intersection of the
finite number of any of its elements is not empty).

{\bf Theorem A.}([3, Theorem 3.3]) \emph{Let $\mathbb{F}$ be a
binary $L$-monad which weakly preserves preimages, and let $X$ be
such that $FX$ is an openly generated (connected) compactum. Then
each map $f:FX\to Y$ with $\mathbb{F}$-convex fibers is 0-soft
(soft) provided $f$ is open.}

By $\exp X$, for any compact $X$, we denote the space of all
nonempty closed subsets of $X$ equipped with the Vietoris topology
(see, e.g., [10]).

In what follows we shall need the characterization of $OS(X)$, given
in [5]. In particular, the following facts take place:

\begin{itemize}
  \item For any $A\in \exp P(X)$ the functional $\nu_A$ given by
$\nu_A(\varphi)=\sup\lbrace \mu(\varphi) \vert \mu\in A\rbrace$,
where $\varphi\in C(X)$, exists and belongs to $OS(X)$. Also
$\nu_A=\nu_{\mathrm{conv}(A)}$ for any $A\in \exp P(X)$ (Proposition
3.2);
  \item Any $\nu\in OS(X)$ coincides with a functional of the form
$\nu_A$, where $A=\lbrace \mu\in P(X)\vert
\mu(\varphi)\le\nu(\varphi) \ \forall \varphi\in C(X)\rbrace$ is a
convex compactum in $P(X)$, in addition, for each $\varphi\in C(X)$
there is $\mu\in A$ such that $\mu(\varphi)=\nu(\varphi)$ (Theorem
3.3);
  \item The correspondence between functionals from $OS(X)$ and
closed convex subsets of $P(X)$ is one-to-one (Theorem 3.4);
  \item For any $f:X\to Y$ and $\nu_A\in OS(X)$ we have
$OS(f)(\nu_A)=\nu_{P(f)(A)}$.
\end{itemize}

\vspace*{\parindent}

{\bf 2. When $OS(X)$ and $OH(X)$ are absolute retracts?}

For any subset $A\subset OH(X)$, we see that $\sup A, \ \inf A$ also
belong to $OH(X)$. Thus, $OH(X)$ is a compact sublattice of
$\prod_{\varphi\in CX}[\min\varphi,\max\varphi]$.

The following statement can be obtained by applying the same
arguments as in [4, Theorem 1].

\vspace*{\parindent}

{\bf Proposition 1.} \emph{For any surjective function $f:X\to Y$
the mapping $OH(f)$ is open if and only if $f$ is open.}

\vspace*{\parindent}

From the remarks on $OS$ made in the first section one can see that
$OS$ is in fact isomorphic to the composition of the functors $cc$
and $P$. Some properties of the functor $cc$ were studied in [8].
For any convex compact $X$, $ccX$ is defined to be the set of all
nonempty closed convex subsets of $X$, $ccX$ is considered as the
subspace of $\exp X$. For any affine mapping $f:X\to Y$ function
$cc(f)$ is given by $cc(f)(A)=f(A)$ where $A\in ccX$. From [8,
Proposition 3.1] and openness of the functor of probability measures
follows

\vspace*{\parindent}

{\bf Proposition 2.} \emph{The functor $OS$ is open, i.e. for any
open mapping $f:X\to Y$ the map $OS(f)$ is open}.

\vspace*{\parindent}

It was shown in [3] that the monad $\mathbb{O}$ generated by the
functor of weakly additive functionals weakly preserves preimages
(Theorem 4.2). Since $\mathbb{OH}$ and $\mathbb{OS}$ are submonads
of $\mathbb{O}$, they weakly preserve preimages as well.

Recall that the notation $\mathbb{L}$ stands for the superextension
monad generated by the superextension functor $\lambda$ (see [10]
for details). For any compact $X$, the space $\lambda X$ has a
functional representation which can be defined by the embedding
$iX:\lambda X \to\prod_{\varphi\in CX}[\min\varphi,\max\varphi]$
such that $iX({\cal A})(\varphi)=\sup\lbrace \inf \varphi (A)\vert
A\in {\cal A}\rbrace$, where ${\cal A}$ is from $\lambda X$ and
$\varphi\in C(X)$. It is easy to see that the image $iX(\lambda X)$
lies in $OH(X)$. Actually, the natural transformation $i=\lbrace
iX\rbrace$ is a monad morphism which embeds the superextension monad
in $\mathbb{OH}$. Therefore, by [3, Theorem 3.2], $\mathbb{OH}$ is
binary.

Now take any openly generated compactum $X$. Whereas the functor
$OH$ is open and the space $OH(X)$ is convex, $OH(X)$ is an openly
generated continuum. From Theorem A we see that whenever
$\mathbb{F}$ is a binary $L$-monad that weakly preserves preimages,
then $F(X)\in AR$ for some compact $X$ provided $F(X)$ is an openly
generated connected compactum. Applying this fact in our case we see
that $OH(X)\in AR$.

Conversely, if we suppose that $OH(X)\in AR$ for some compact $X$,
then an argumentation similar to that of [4, Theorem 2] provides
that $X$ is an openly generated compactum.

We therefore obtain the following fact:

\vspace*{\parindent}

{\bf Theorem 1.} \emph{$OH(X)$ is an absolute retract if and only if
$X$ is an openly generated compactum}.

So what we get is that $OH(X)$ can be an $AR$ even when the weight
of $X$ exceeds $\omega_1$. The same could be said on some other
functors which generate $L$-monads and contain $\mathbb{L}$ as
submonad, for instance $G,O$, $\lambda$ by itself. The functor $OS$
seems to be closer to $P$. It does not give an $AR$ in weights
higher than $\omega_1$:

\vspace*{\parindent}

{\bf Proposition 3.} \emph{$OS(X)$ is an absolute retract if and
only if $X$ is openly generated with {\it w}($X$)$\le \omega_1$}.

\emph{Proof.} Follows from the results of [8], namely [8, Theorem
4.1] combined with results of [7] providing that a statement
analogous to that of the proposition holds for the functors $cc$ and
$P$.

\vspace*{\parindent}

{\bf Corollary 1.} \emph{There is no monad embedding
$i:\mathbb{L}\hookrightarrow\mathbb{OS}$.}

Indeed, assuming the contrary, we would obtain that $\mathbb{OS}$ is
binary. Therefore, according to [3, Theorem 3.3], the space
$OS(D^{\omega_2})$, for example, must be an absolute retract, a
contradiction.

\vspace*{\parindent}

{\bf 3. The softness of multiplication maps for $\mathbb{OH}$ and
$\mathbb{OS}$}.

\vspace*{\parindent}

{\bf Theorem 2.} \emph{If the multiplication map $\mu_{OS} X$ for
$\mathbb{OS}$ is soft then $X$ is metrizable}.

\emph{Proof}. Suppose that $X$ is not metrizable and $\mu_{OS} X$ is
soft. Use [9, Theorem 3] to obtain that $X$ is openly generated.

Represent $X$ as the limit of an $\omega$-system $S=\lbrace
X_\alpha, p_{\alpha}^{\beta}, \cal{A}\rbrace$ with open bonding
maps. Whereas $\mu X$ is soft, we can assume that all limit diagrams

\begin{displaymath}
\xymatrix{
  OS^2(X) \ar[r]^{OS^2(p_\alpha)} \ar[d]_{\mu X} & OS^(X_\alpha)
   \ar[d]^{\mu X_\alpha}\\
  OS(X) \ar[r]^{OS(p_\alpha)} & OS(X_\alpha)}
\end{displaymath}
are soft [9, Theorem 2], hence open.

Now our aim is to obtain $\alpha_0\in \cal{A}$ and an accumulation
point $x\in X_{\alpha_0}$ such that $p_{\alpha_0}^{-1}(x)$ contains
more than one point. The weight of $X$ is uncountable, so its
character is uncountable too, since $w(X)=\chi(X)$ for any openly
generated compactum [4]. Choose $x_0\in X$ with $\chi(x_0,X)>\omega$
and some $\alpha\in\cal{A}$, put $x_\alpha=p_\alpha(x_0)$. Then
$p_{\alpha}^{-1}(x_\alpha)$ contains more than one point, otherwise
$x_0$ would have the countable character. If $x_\alpha$ is not
isolated, Then $x_\alpha$ is the required point. Suppose that
$x_\alpha$ is isolated. Consider $x_1\in p_{\alpha}^{-1}(x_\alpha)$
distinct from $x_0$. We can choose $\alpha_1> \alpha$ with
$p_{\alpha_1}(x_1)\neq p_{\alpha_1}(x_0)$. Again
$p_{\alpha_1}^{-1}(x_0)$ is not a singleton, and if
$p_{\alpha_1}(x_0)$ is an accumulation point, we are done. Assume
the opposite. Take any $x_2\in p_{\alpha_1}^{-1}(p_{\alpha_1}(x_0))$
with $x_2\neq x_0$ and $\alpha_2> \alpha_1$ such that
$p_{\alpha_2}(x_2)\neq p_{\alpha_2}(x_0)$ and continue the process
as described above. If on any step $i$ the point $p_{\alpha_i}(x_0)$
is not an accumulation point, we obtain the sequence $\lbrace
x_i\rbrace_{i\in\mathbb{N}}$ of points in $X$ and the up-directed
chain of elements $\lbrace \alpha_i\rbrace_{i\in\mathbb{N}}$ of
$\cal{A}$ which has the least upper bound $\alpha_0\in \cal{A}$.
Then the space $X_{\alpha_0}$ is the limit of the inverse system
$\lbrace X_{\alpha_i},p_{\alpha_i}^{\alpha_j},i\le j\rbrace$ and
$\lim_{i\to\infty}p_{\alpha_0}(x_i)=p_{\alpha_0}(x_0)$. Indeed, the
family $\lbrace
(p_{\alpha_i}^{\alpha_0})^{-1}(p_{\alpha_i}(x_0))\rbrace$ forms a
base of neighborhoods at $p_\alpha(x_0)$, and for any such
$(p_{\alpha_i}^{\alpha_0})^{-1}(p_{\alpha_i}(x_0))$ we see that
$p_{\alpha_0}(x_j)$ is contained in it for all $j\ge i$. Therefore,
$\alpha_0\in\cal{A}$ and $x = p_{\alpha_0}(x_0)$ chosen above are as
required.

According to our assumption, the diagram

\begin{displaymath}
\xymatrix{
  OS^2(X) \ar[r]^{OS^2(p_{\alpha_0})} \ar[d]_{\mu_{OS} X} & OS^(X_{\alpha_0})
   \ar[d]^{\mu_{OS} X_{\alpha_0}}\\
  OS(X) \ar[r]^{OS(p_{\alpha_0})} & OS(X_{\alpha_0})}
\end{displaymath}
is open.

Consider the accumulation point $x\in X_{\alpha_0}$ chosen above and
distinct $y_1,y_2\in p_{\alpha_0}^{-1}(x)$. Let $\lbrace
x_i\rbrace_{\i\in I}$ be the net converging to $x$. Choose $y_i\in
p_{\alpha_0}^{-1}(x_i)$ for every $i\in I$ the way that
$\{y_i\}_{i\in I}$ would converge to $x_1$.

Denote $\nu =(\delta_{y_1}+\delta_{y_2})/ 2$ and ${\cal V} =
\Delta_{\delta_x}$. Then the net $(\nu_i,{\cal V}_i)$=
$((\delta_{y_i}+\delta_{y_2})/ 2,\Delta_{(\delta_x+\delta_{x_i})/
2})$ converges to $(\nu,{\cal V})$. To obtain a contradiction with
openness of $\chi$, and therefore softness of $\mu_{OS} X$, show
that the inverse of the characteristic map of the considered diagram
is not continuous. Indeed, let us consider
$\chi^{-1}(((\delta_{y_i}+
\delta_{y_2})/2,\Delta_{(\delta_x+\delta_{x_i})/
2})=(OS^2(p_{\alpha_0}))^{-1}({\cal V}_i)\cap (\mu_{OS}
X)^{-1}(\nu_i)$. Suppose that $\Theta\in
(OS^2(p_{\alpha_0}))^{-1}({\cal V}_i)$. Then $\mathrm{supp} \
\Theta$ is in $OS((p_{\alpha_0})^{-1}(x)\cup
(p_{\alpha_0})^{-1}(x_i))$. Now take any $\Theta =
\Theta_A\in(\mu_{OS} X)^{-1}(\nu_i)$ (recall that any functional
$\eta\in OS(Y)$ is of the form $\eta = \eta_B$, where $B\in
ccP(Y)$). We want to show that $\mathrm{supp} \ \Theta_A$ is in
$OS(\{y_i,y_2\})$. Indeed, assuming the contrary, we obtain that
there is a measure $M\in A$ that is not supported on
$OS(\{y_i,y_2\})$, therefore there exists $\theta\notin
OS(\{y_i,y_2\})$ from the support of $M$. So we can choose a
function $\varphi\in C(X)$ which is zero at $\{y_2,y_i\}$,
$\theta(\varphi)>0$ and $0_X\le \varphi$. Since $\pi_\varphi$ is
continuous, there exists a closed neighborhood $V$ of $\theta$ on
which $\pi_\varphi$ is strictly greater than zero. Also,
$\mathrm{supp} \ M\cap V\neq \O$, so $M(V)>0$. This implies
$M(\pi_\varphi)>0$, and hence $\mu_{OS}
X(\Theta_A)(\varphi)=\Theta_A(\pi_\varphi)=\sup \{ M(\pi_\varphi) \
\vert \ M\in A \}>0$, whereas $\nu_i(\varphi)=0$ which gives us
$\mu_{OS} X(\Theta_A)\neq \nu_i$. That's why any $\Theta\in
(OS^2(p_{\alpha_0}))^{-1}({\cal V}_i)\cap (\mu_{OS} X)^{-1}(\nu_i)$
must be supported on $OS(\{y_i,y_2\})$. The only such functional
$\Theta$ which also satisfies the condition
$\chi(\Theta)=((\delta_{y_i}+\delta_{y_2})/2,\Delta_{(\delta_x+\delta_{x_i})/
2})$ is the measure $\Delta_{(\delta_{y_i}+\delta_{y_2}/2}$.
Therefore, there is some neighborhood $V_1$ of the functional
$(\Delta_{\delta_{y_1}}+\Delta_{\delta_{y_2}})/
2\in\chi^{-1}(\nu,{\cal V})$ that contains no elements of the form
$\Delta_{(\delta_{y_i}+\delta_{y_2})/ 2}$ starting from some $i_0\in
I$, hence $\chi^{-1}$ is not continuous, and the diagram is not
open, a contradiction with the initial assumption. Theorem is
proved.

\vspace*{\parindent}

The following are the results for $\mathbb{OH}$ which show that it
behaves the same way as the monad $\mathbb{O}$.

\vspace*{\parindent}

{\bf Theorem 3.} \emph{$\mu_{OH} X$ is open for any compactum $X$}.

Proof of Theorem 3 is the same as that for $\mathbb{O}$ [3].

{\bf Theorem 4.} \emph{$\mu_{OH} X$ is soft if and only if $X$ is an
openly generated compactum}.

\emph{Proof}. \emph{Necessity}. Let $X=\lim {\cal S}$, where ${\cal
S}=\{X_\alpha,p_\alpha,{\cal A}\}$ is an $\omega$-system consisting
of metrizable compacta and epimorphisms. The mapping $\mu_{OH} X$ is
soft, hence we can assume that all the limit diagrams of the form

\begin{displaymath}
\xymatrix{
  OH^2(X) \ar[r]^{OH^2(p_{\alpha})} \ar[d]_{\mu_{OH} X} & OH^(X_{\alpha})
   \ar[d]^{\mu_{OH} X_{\alpha}}\\
  OH(X) \ar[r]^{OH(p_{\alpha})} & OH(X_{\alpha})}
\end{displaymath}
are open. Assume that $X$ is not openly generated, so that there
exists $\alpha\in {\cal A}$ such that $p_\alpha$ is not open. Then
by Proposition 2 the mapping $OH(p_\alpha)$ is not open. Therefore,
there is a functional $\nu\in OH(X_\alpha)$ and a net
$\{\nu_i\}_{i\in I}$ converging to $\nu$ such that the net
$OH(p_\alpha)^{-1}(\nu_i)$ converges to some $A\neq
OH(p_\alpha)^{-1}(\nu)$. We have that $A\subset
OH(p_\alpha)^{-1}(\nu)$. Choose two comparable elements $\theta_1\in
A$ and $\theta_2\in OH(p_\alpha)^{-1}(\nu)\backslash A$. Let, for
example, $\theta_1\le \theta_2$. Let $\{\theta_i\}$ be a net
converging to $\theta_2$ such that $\theta_i\in
OH(p_\alpha)^{-1}(\nu_i)$ for all $i\in I$. We see that the net
$\{(\theta_i,\eta OH(X_\alpha)(\nu_i))\}$ converges to $(\theta_2,
\eta OH(X_\alpha)(\nu))$. Now let ${\cal V}\in OH^2(X)$ be a
functional such that ${\cal V}(\Phi)=\max
\{\Phi(\theta_1),\Phi(\theta_2)\}$. Then $\chi({\cal
V})=(\theta_2,\eta OH(X_\alpha)(\nu))$, where $\chi$ is the
characteristic map of the diagram. Choose $\Phi\in C(OH(X))$ with
$\Phi(\theta_2)=1$ and $\Phi(\theta)=0$ for any $\theta\in A$. Then
we may take that $\Phi(\theta)\le\1/2$ for any
$OH(p_\alpha)^{-1}(\nu_i)$, hence, using [4, Lemma 2], we get that
$\Theta(\Phi)\le 1/2$ for all $\Theta\in (OH^2(p_\alpha))^{-1}(\eta
OH(X_\alpha)(\nu_i))$. Thus we obtained an open neighborhood of
${\cal V}$ of the form $V=\{\Theta\in OH^2(X) \ \ \vert
\Theta(\Phi)>1/2\}$ with $V\cap \chi^{-1}(\theta_i,\eta
OH(X_\alpha)(\nu_i))=\O$, a contradiction which shows that $X$ must
be openly generated.

\emph{Sufficiency}. The monad $\mathbb{OH}$ is a binary monad which
weakly preserves preimages. Since $\mu_{OH} X:OH^2(X)\to OH(X)$ is
an open $\mathbb{OH}$-algebras morphism  and $OH^2(X)$ is openly
generated (by Theorem 3), the softness of $\mu_{OH} X$ follows from
Theorem A. The statement is proved.

 \linespread{1} \small{

}

\end{document}